\documentclass[11pt]{article}
\renewcommand{\baselinestretch}{1.2}
\newcommand{\dated}{\mbox{} \hfill {\small [{\tt \today}]}} 
\usepackage{amsthm,enumerate}
\theoremstyle{plain}
\newtheorem{theorem}{Theorem}[section]
\newtheorem{lemma}[theorem]{Lemma}
\newtheorem{corollary}[theorem]{Corollary}
\newtheorem{proposition}[theorem]{Proposition}
\theoremstyle{definition}
\newtheorem{definition}[theorem]{Definition}
\theoremstyle{remark}
\newtheorem*{remark}{Remark}
\newtheorem*{example}{Example}
\newtheorem*{rems}{Remarks}
\newtheorem*{exs}{Examples}
\newenvironment{remarks}{\begin{rems}\begin{enumerate}}{\end{enumerate}\end{rems}}

\newenvironment{items}{\begin{enumerate}[\rm (i)]}{\end{enumerate}}

 \usepackage{amsmath,amssymb,amsfonts,diagrams}
\newenvironment{keywords}{\noindent\small {\it Keywords\/}:}{\vskip 4pt}
\newenvironment{classification}{\noindent\small 2000 {\it Mathematics Subject
Classification\/}:}{\vskip 12pt}

\renewcommand{\iff}{\quad\Longleftrightarrow\quad}

\newcommand{\comps}{{\mathbb C}}

\newcommand{\free}{{\mathbb F}}

\newcommand{\tensor}{\otimes}

\newcommand{\id}{{\mathrm{id}}}
\newcommand{\cb}{{\mathrm{cb}}}

\newcommand{\A}{{\mathfrak A}}

\newcommand{\Hilbert}{{\mathfrak H}}

\theoremstyle{plain}
\newtheorem*{gilbert}{Gilbert's Theorem}
\newcommand{\Mcb}{\mathit{Mcb}}
\newcommand{\MA}{\mathit{MA}}
\title{Norm one idempotent $cb$-multipliers \\ with applications to the Fourier algebra \\
in the $cb$-multiplier norm}
\author{\textit{Brian E.\ Forrest}\thanks{Research supported by NSERC.} \and \textit{Volker Runde}\thanks{Research supported by NSERC.}}
\date{}
\begin{document}
\maketitle
\begin{abstract}
For a locally compact group $G$, let $A(G)$ be its Fourier algebra, let $M_{cb}A(G)$ denote the completely bounded multipliers of $A(G)$, and let $A_{\Mcb}(G)$ stand for the closure of $A(G)$ in $M_{cb}A(G)$. We characterize the norm one idempotents in $M_{cb}A(G)$: the indicator function of a set $E \subset G$ is a norm one idempotent in $M_{cb}A(G)$ if and only if $E$ is a coset of an open subgroup of $G$. As applications, we describe the closed ideals of $A_{\Mcb}(G)$ with an approximate identity bounded by $1$, and we characterize those $G$ for which $A_{\Mcb}(G)$ is $1$-amenable in the sense of B.\ E.\ Johnson. (We can even slightly relax the norm bounds.)
\end{abstract}
\section*{Introduction}
The Fourier algebra $A(G)$ and Fourier--Stieltjes algebra $B(G)$ of a locally compact group $G$ were introduced by P.\ Eymard in \cite{Eym}. If $G$ is abelian with dual group $\hat{G}$, these algebras are isometrically isomorphic to $L^1(\hat{G})$, the group algebra of $\hat{G}$, and $M(\hat{G})$, the measure algebra of $\hat{G}$, via the Fourier and Fourier--Stieltjes transform, respectively. For abelian $G$, the idempotent elements in $B(G) \cong M(\hat{G})$ were described by P.\ J.\ Cohen (\cite{Coh}): the indicator function $\chi_E$ of $E \subset G$ lies in $B(G)$ if and only if $E$ belongs to the \emph{coset ring} $\Omega(G)$ of $G$, i.e., the ring of sets generated by the cosets of the open subgroups of $G$. Later, B.\ Host showed that this characterization of the idempotents in $B(G)$ holds true for general locally compact groups $G$ (\cite{Hos}).
\par 
In \cite{FKLS}, the Cohen--Host idempotent theorem was crucial in characterizing, for amenable $G$, those closed ideals of $A(G)$ that have a bounded approximate identity, and in \cite{FR} and \cite{Run}, the authors made use of it to characterize those $G$ for which $A(G)$ is amenable in the sense of B.\ E.\ Johnson (\cite{Joh}).
\par 
Besides the given norm on $A(G)$, there are other---from certain points of view even more natural---norms on $A(G)$. Recall that a \emph{multiplier} of $A(G)$ is a function $\phi$ on $G$ with $\phi A(G) \subset A(G)$. It is immediate from the closed graph theorem that each multiplier $\phi$ of $A(G)$ induces a bounded multiplication operator $M_\phi$ on $A(G)$; the operator norm on the multipliers turns them into a Banach algebra. Trivially, $A(G)$ embeds contractively into its multipliers, but the multiplier norm on $A(G)$ is equivalent to the given norm if and only if $G$ is amenable (\cite{Los}).
\par 
An even more natural norm on $A(G)$ arises if we take into account that $A(G)$, being the predual of a von Neumann algebra, has a canonical operator space structure. (Our default reference for operator spaces is \cite{ER}.) This makes it possible to consider the \emph{completely bounded multipliers}---\emph{$cb$-multipliers} in short---of $A(G)$ as
\[
  M_{cb}A(G) := \{ \phi : G \to \comps : \text{$M_\phi \!: A(G) \to A(G)$ is completely bounded} \}.
\]
For $\phi \in M_{cb}A(G)$, we denote the completely bounded operator norm of $M_\phi$ by $\| \phi \|_{\Mcb}$. It is not difficult to see that $B(G)$ embeds completely contractively into $M_{cb}A(G)$. However, equality holds if and only if $G$ is amenable. In fact, $G$ is amenable if and only if $\| \cdot \|_{Mcb}$ and the given norm on $A(G)$ are equivalent. (For a discussion of these facts with references to the original literature, see \cite{Spr}.)
\par 
Let $A_{\Mcb}(G)$ denote the closure of $A(G)$ in $M_{cb}A(G)$ (see \cite{For2} for some properties of this algebra). For certain non-amenable $G$, the (completely contractive) Banach algebra $A_{\Mcb}(G)$ is better behaved than $A(G)$. For instance, $A(G)$ has a bounded approximate identity if and only if $G$ is amenable (\cite{Lep}); in particular, if $G$ is $\free_2$, the free group in two generators, then $A(G)$ is not operator amenable. On the other hand, $A_{\Mcb}(\free_2)$ has a bounded approximate identity (\cite{dCH}) and even is operator amenable (\cite{FRS}) in the sense of \cite{Rua}.
\par 
Juxtaposing the main results of \cite{FR} and \cite{FRS}, the question arises immediately whether $A_{\Mcb}(\free_2)$ is amenable in the classical sense of \cite{Joh}, and it is this question that has motivated the present note. The proof of the main result of \cite{FR}--as well as its alternative proof in \cite{Run}---rests on the Cohen--Host idempotent theorem. Attempting to emulated these proofs with $A_{\Mcb}(G)$ in place of $A(G)$ leads to the problem whether certain idempotent functions can lie in $M_{cb}A(G)$. The main problem is that the Cohen--Host theorem is no longer true with $M_{cb}A(G)$ replacing $B(G)$: as M.\ Leinert showed (\cite{Lei}), there are sets $E \subset \free_2$ such that $\chi_E \in M_{cb}A(G) \setminus B(G)$. 
\par
The main result of this note is that, even though $M_{cb}(G)$ may have more idempotents than $B(G)$, both algebras do have the same \emph{norm one} idempotents. With this result, we can then characterize the closed ideals in $A_{\Mcb}(G)$ having an approximate identity bounded by one as well as those $G$ for which $A_{\Mcb}(G)$ is $1$-amenable. (Due to the useful fact that idempotent Schur multipliers of norm less than $\frac{2}{\sqrt{3}}$ must have norm one, we can even work with slightly relaxed norm bounds.)
\section{The norm one idempotents of $M_{cb} A(G)$}
For a locally compact group $G$, the functions in $B(G)$ can be described as coefficient functions of unitary representations of $G$ (see \cite{Eym}). A related characterization, which immediately yields the contractive inclusion $B(G) \subset M_{cb}A(G)$, is the following theorem due to J.\ Gilbert (\cite{Gil}; for a more accessible proof, see \cite{Jol}):
\begin{gilbert}
Let $G$ be a locally compact group. Then the following are equivalent for $\phi \!: G \to \comps$:
\begin{items}
\item $\phi \in M_{cb}A(G)$;
\item there are are a Hilbert space $\Hilbert$ and bounded, continuous functions $\boldsymbol{\xi}, \boldsymbol{\eta} \!: G \to \Hilbert$ such that
\begin{equation} \label{gileq}
  \phi(xy^{-1}) = \langle \boldsymbol{\xi}(x), \boldsymbol{\eta}(y) \rangle \qquad (x,y \in G).
\end{equation}
\end{items}
Moreover, if $\phi \in M_{cb}A(G)$ and $\boldsymbol{\xi}$ and $\boldsymbol{\eta}$ are Hilbert space valued, bounded, continuous functions on $G$ satisfying \emph{(\ref{gileq})} then 
\begin{equation} \label{gilineq}
  \| \phi \|_{\Mcb} \leq \| \boldsymbol{\xi} \|_\infty \| \boldsymbol{\eta} \|_\infty
\end{equation} 
holds, and $\boldsymbol{\xi}$ and $\boldsymbol{\eta}$ can be chosen such that we have equality in \emph{(\ref{gilineq})}.
\end{gilbert}
\par 
The following extends \cite[Theorem 2.1]{IS}:
\begin{theorem} \label{idemthm}
Let $G$ be a locally compact group. Then the following are equivalent for $E \subset G$:
\begin{items}
\item $\chi_E \in B(G)$ with $\| \chi_E \|_{B(G)} = 1$; 
\item $\chi_E \in M_{cb}A(G)$ with $\| \chi_E \|_{\Mcb} = 1$;
\item $E$ is a coset of an open subgroup.
\end{items}
\end{theorem}
\begin{proof}
(i) $\Longrightarrow$ (ii) is clear, and (iii) $\Longrightarrow$ (i) is the easy part of \cite[Theorem 2.1]{IS}.
\par
(ii) $\Longrightarrow$ (iii): Obviously, $E$ is open. If $x \in E$, then $x^{-1}E$ contains $e$ and satisfies  $\| \chi_{x^{-1}E} \|_{\Mcb} =1$. Hence, we can suppose without loss of generality that $e \in E$: otherwise, replace $E$ by $x^{-1}E$ for some $x \in E$. We shall show that $E$ is a subgroup of $G$.
\par
By Gilbert's Theorem, there are a Hilbert space $\Hilbert$ and bounded, continuous functions $\boldsymbol{\xi}, \boldsymbol{\eta} \!: G \to \Hilbert$ with $1 = \| \boldsymbol{\xi} \|_\infty \| \boldsymbol{\eta} \|_\infty$ such that
\begin{equation} \label{Gilbert}
  \chi_E(xy^{-1}) = \langle \boldsymbol{\xi}(x), \boldsymbol{\eta}(y) \rangle \qquad (x,y \in G);
\end{equation}
of course, we can suppose that both $\| \boldsymbol{\xi} \|_\infty = \| \boldsymbol{\eta} \|_\infty =1$. In view of (\ref{Gilbert}) and the Cauchy--Schwarz inequality, we obtain
\[
  xy^{-1} \in E \iff \langle \boldsymbol{\xi}(x), \boldsymbol{\eta}(y) \rangle = 1 \iff \boldsymbol{\xi}(x) = \boldsymbol{\eta}(y)
  \qquad (x,y \in G).
\]
As $e \in E$, this means, in particular, that $\boldsymbol{\xi}(e) = \boldsymbol{\eta}(e) =:\xi$, so that
\[
  E = \{ x \in G : \boldsymbol{\xi}(x) = \xi \} = \{ y \in G : \boldsymbol{\eta}(y^{-1}) = \xi \}.
\]
Hence, if $x,y \in E$, we get
\[
  \chi_E(xy) = \langle \boldsymbol{\xi}(x), \boldsymbol{\eta}(y^{-1}) \rangle = \langle \xi, \xi \rangle = 1, 
\]
so that $xy \in E$. Consequently, $E$ is a subsemigroup of $G$.
\par 
Let $x \in E$. Applying the preceding argument to $x^{-1}E$ instead of $E$, we see that $x^{-1}$ is a subsemigroup of $G$; since $e \in E$, we have, in particular, $x^{-1} x^{-1} \in x^{-1}E$, which means that $x^{-1} \in E$. 
\par 
All in all, $E$ is a subgroup of $G$.
\end{proof}
\begin{remark}
Let $\MA(G)$ denote the algebra of all multipliers of $A(G)$. Defining $\| \phi \|_M$ as the operator norm of $M_\phi$, we obtain a Banach algebra norm on $\MA(G)$; obviously, $M_{cb}A(G)$ embeds contractively into $\MA(G)$. Hence, every norm one idempotent in $M_{cb}A(G)$ is a norm one idempotent in $\MA(G)$. By \cite{Boz}, $M_{cb}A(\free_2) \subsetneq \MA(\free_2)$ holds, and---as M.\ Bo\.zejko communicated to the second named author---, there are sets $E \subset \free_2$ such that $\chi_E \subset \MA(\free_2) \setminus M_{cb}A(\free_2)$. We do not know if such $E$ can be chosen such that $\| \chi_E \|_M = 1$.
\end{remark}
\par 
By \cite{BF}, the elements of $M_{cb}A(G)$ are precisely the so-called \emph{Herz--Schur multipliers} of $A(G)$. For discrete $G$, the powerful theory of \emph{Schur multipliers} (see \cite{Pis} for an account) can thus be applied to the study of $M_{cb}A(G)$. By \cite{Liv} (see also \cite{KP}), an idempotent Schur multiplier of ${\cal B}(\ell^2(\mathbb{I}))$, where $\mathbb I$ is any index set, with norm greater than $1$ must have norm at least $\frac{2}{\sqrt{3}}$. Hence, we obtain:
\begin{corollary} \label{idemcor}
Let $G$ be a group. Then the following are equivalent for $E \subset G$:
\begin{items}
\item $\chi_E \in B(G)$ with $\| \chi_E \|_{B(G)} = 1$; 
\item $\chi_E \in M_{cb}A(G)$ with $\| \chi_E \|_{\Mcb} = 1$;
\item $\chi_E \in M_{cb}A(G)$ with $\| \chi_E \|_{\Mcb} < \frac{2}{\sqrt{3}}$;
\item $E$ is a coset of a subgroup.
\end{items}
\end{corollary}
\section{Ideals of $A_{\Mcb}(G)$ with approximate identities bounded by $C < \frac{2}{\sqrt{3}}$}
Let $G$ be a locally compact group. In \cite{FKLS}, the first named author with E.\ Kaniuth, A.\ T.-M.\ Lau, and N.\ Spronk characterized, for amenable $G$, those closed ideals of $A(G)$ that have bounded approximate identities in terms of their hulls. Previously, he had obtained a similar characterization of those closed ideals of $A(G)$ that have approximate identities bounded by one without any amenability hypothesis for $G$ (\cite[Propositon 3.12]{For1}).
\par
In this section, we use Theorem \ref{idemthm} (or rather Corollary \ref{idemcor}) to prove an analog of \cite[Propositon 3.12]{For1} for $A_{\Mcb}(G)$.
\par 
Let $H$ be an open subgroup of $G$. It is well known that we can isometrically identify $A(H)$ with the closed ideal of $A(G)$ consisting of those functions whose support lies in $H$; with a little extra effort, one sees that this identification is, in fact, a complete isometry (\cite[Proposition 4.3]{FW}). From there, it is not difficult to prove the analogous statement for $A_{\Mcb}(G)$: there is a canonical isometric isomorphism between $A_{\Mcb}(H)$ and those functions in $A_{\Mcb}(G)$ with support in $H$.
\par 
Given a closed ideal $I$ of $A_{\Mcb}(G)$, we define its \emph{hull} to be
\[
  h(I) := \{ x \in G : \text{$f(x) = 0$ for all $f \in I$} \}.
\]
If $E \subset G$ is closed, we set
\[
  I(E) := \{ f \in A_{\Mcb}(G) : \text{$f(x) = 0$ for all $x \in E$} \},
\]
which is a closed ideal of $A_{\Mcb}(G)$ such that $h(I(E)) = E$.
\par 
Since translation by a group element is an isometric algebra automorphism of $A_{\Mcb}(G)$, we have in view of the preceding discussion:
\begin{proposition} \label{isoprop}
Let $G$ be a locally compact group, let $H$ be an open subgroup of $G$, and let $x \in G$. Then we have an isometric algebra isomorphism between $A_{\Mcb}(H)$ and $I(G \setminus x H)$. 
\end{proposition}
\par 
Our main result in this section is the following:
\begin{theorem} \label{baithm}
Let $G$ be a locally compact group. Then the following are equivalent for a closed ideal $I$ of $A_{\Mcb}(G)$ and $C \in \left[1,\frac{2}{\sqrt{3}} \right)$:
\begin{items}
\item $I$ has an approximate identity bounded by $C$;
\item $I = I(G \setminus xH)$, where $x \in G$ and $H$ is an open subgroup of $G$ such that $A_{\Mcb}(H)$ has an approximate identity bounded by $C$.
\end{items}
\end{theorem}
\begin{proof}
(ii) $\Longrightarrow$ (i) is an immediate consequence of Proposition \ref{isoprop}.
\par 
(i) $\Longrightarrow$ (ii): Let $( e_\alpha )_\alpha$ be an approximate identity for $I$ bounded by $C$. By \cite[Corollary 6.3(i)]{Spr}, $M_{cb}A(G)$ embeds (completely) isometrically into $M_{cb}A(G_d)$, where $G_d$ stands for the group $G$ equipped with the discrete topology; we may thus view $( e_\alpha )_\alpha$ as a bounded net in $M_{cb}A(G_d)$. It is easy to see that $( e_\alpha )_\alpha$ converges to $\chi_{G \setminus h(I)}$ pointwise on $G$ and thus in $\sigma(\ell^\infty(G),\ell^1(G))$. With the help of \cite[Lemma 1.9]{dCH}, we conclude that $\chi_{G \setminus h(I)} \in M_{cb}A(G_d)$ with $\| \chi_{G \setminus h(I)} \|_{\Mcb} \leq C$; 
hence, $\chi_{G \setminus h(I)}$ is an idempotent in $M_{cb}A(G_d)$ of norm strictly less than $\frac{2}{\sqrt{3}}$. By Corollary \ref{idemcor}, this means that $G \setminus h(I)$ is of the form $xH$ for $x \in G$ and a subgroup $H$ of $G$ and thus $h(I) = G \setminus xH$. Since $h(I)$ is closed, $xH$---and thus $H$---must be open. By \cite[Proposition 2.2]{FRS}, the Banach algebra $A_{\Mcb}(H)$ is Tauberian. By Proposition \ref{isoprop}, this means that the set $G \setminus xH$ is of synthesis for $A_{\Mcb}(G)$, so that $I = I(G \setminus xH)$. Finally, Proposition \ref{isoprop} again yields that $A_{\Mcb}(H)$ has an approximate identity bounded by $C$.
\end{proof}
\par 
In \cite{CH}, locally compact groups $G$ such that $A(G)$ has an approximate identity bounded in $\| \cdot \|_{\Mcb}$ were called \emph{weakly amenable}; this is equivalent to $A_{\Mcb}(G)$ having an approximate identity (\cite[Proposition 1]{For2}). For instance, $\free_2$ is weakly amenable \cite[Corollary 3.9]{dCH} without being amenable. Both \cite[Corollary 3.9]{dCH} and \cite[Theorem 2.7]{FRS} suggest that, for weakly amenable, but not amenable $G$, the Banach algebra $A_{\Mcb}(G)$ is a more promising object of study than $A(G)$. In view of \cite[Proposition 3.13]{For1} and Theorem \ref{baithm}, one is thus tempted to ask whether a suitable version of \cite[Theorem 2.3]{FKLS} holds for $A_{\Mcb}(G)$ and weakly amenable $G$: a closed ideal $I$ of $A_{\Mcb}(G)$ has a bounded approximate identity if and only $I = I(E)$ for some closed $E \in \Omega(G_d)$. 
\par 
We conclude this section with an example which shows that the characterization of the closed ideals of $A_{\Mcb}(G)$ with a bounded approximate identity for weakly amenable, but not amenable $G$ cannot be as elegant as for amenable $G$:
\begin{example}
Let $E \subset \free_2$ be such that $\chi_E \in M_{cb}A(\free_2)$, but $E \notin \Omega(\free_2)$: such $E$ exists by \cite{Lei}. Let $I = I(E)$. Then $I = (1-\chi_E)A_{\Mcb}(\free_2)$ is completely complemented in $A_{\Mcb}(\free_2)$. Since $A_{\Mcb}(\free_2)$ is operator amenable by \cite[Theorem 2.7]{FRS}, it follows from \cite[Theorem 2.3.7]{LoA}---with operator space overtones added---that $I$ has a bounded approximate identity even though $h(I) = E \notin \Omega(\free_2)$. 
\end{example}
\section{Amenability of $A_{\Mcb}(G)$}
Recall the definition of an amenable Banach algebra.
\par 
Given a Banach algebra $\A$, let $\A \tensor^\gamma \A$ denote the Banach space tensor product of $\A$ with itself. The Banach space $\A \tensor^\gamma \A$ becomes a Banach $\A$-bimodule via
\[
  a \cdot (x \tensor y) := ax \tensor y \quad\text{and}\quad (x \tensor y) \cdot a := x \tensor ya
  \qquad (a,x,y \in \A).
\]
Let $\Delta \!: \A \tensor^\gamma \A \to \A$ denote the bounded linear map induced by multiplication, i.e., $\Delta(a \tensor b) =ab$ for $a,b \in \A$.
\begin{definition} \label{amdef}
A Banach algebra $\A$ is called \emph{$C$-amenable} with $C \geq 1$ if it has an \emph{approximate diagonal} bounded by $C$, i.e., a net $( \boldsymbol{d}_\alpha )_\alpha$ in $\A \tensor^\gamma \A$ bounded by $C$ such that
\[
  a \cdot \boldsymbol{d}_\alpha - \boldsymbol{d}_\alpha \cdot a \to 0 \qquad (a \in \A)
\]
and
\begin{equation} \label{diag2}
  a \Delta \boldsymbol{d}_\alpha \to a \qquad (a \in \A).
\end{equation}
We say that $\A$ is \emph{amenable} if there is $C \geq 1$ such that $\A$ is $C$-amenable.
\end{definition}
\begin{remark}
Remark is not the original definition of an amenable Banach algebra from \cite{Joh}, but equivalent to it (\cite{Joh2}). The idea of considering bounds for approximate diagonals seems to originate in \cite{Joh3}.
\end{remark}
\par 
The question as to which locally compact groups $G$ have an amenable Fourier algebra was first studied in depth in \cite{Joh3}. Until then, it was widely believed---probably with an eye on \cite{Lep}---that these $G$ were precisely the amenable ones. In \cite{Joh3}, however, Johnson exhibited compact groups $G$---such as $\mathrm{SO}(3)$---for which $A(G)$ is not amenable. Eventually, the authors showed that $A(G)$ is amenable if and only if $G$ is almost abelian, i.e., has an abelian subgroup of finite index (\cite[Theorem 2.3]{FR}; see also \cite{Run}).
\par 
A crucial r\^ole in the proofs in both \cite{FR} and \cite{Run} is played by the \emph{anti-diagonal} of $G$; it is defined as
\[
  \Gamma := \{ (x,x^{-1}) : x \in G \}.
\]
Its indicator function $\chi_\Gamma$ lies $B(G_d \times G_d)$ if and only if $G$ is almost abelian (\cite[Proposition 3.2]{Run}). If $G$ is locally compact such that $A(G)$ is amenable, then $\chi_\Gamma$ lies in $B(G_d \times G_d)$ (\cite[Lemma 3.1]{Run}), forcing $G$ to be almost abelian.
\par 
For any $f \!: G \to \comps$, we define $\check{f} \!: G \to \comps$ by letting $\check{f}(x) := f(x^{-1})$. We denote the map assigning $\check{f}$ to $f$ by $\mbox{\ }^\vee$; it is an isometry on $A(G)$, but completely bounded if and only if $G$ is almost abelian (\cite[Proposition 1.5]{FR}): this fact is crucial for characterizing those $G$ with an amenable Fourier algebra as the almost abelian ones (see both \cite{FR} and \cite{Run}).
\par 
Since $\mbox{\ }^\vee$ need not be completely bounded, it is not obvious that $\mbox{\ }^\vee$ is an isometry---or even well defined---on $A_{\Mcb}(G)$. Nevertheless, these are true:
\begin{lemma} \label{checklem}
Let $G$ be a locally compact group. Then $\mbox{\ }^\vee$ is an isometry on $M_{cb}A(G))$ leaving $A_{\Mcb}(G)$ invariant.
\end{lemma}
\begin{proof}
Since $\mbox{\ }^\vee$ leaves $A(G)$ invariant, it is clear that it leaves $A_{\Mcb}(G)$ invariant once we have established that it is isometric on $M_{cb}A(G)$.
\par 
Let $\phi \in M_{cb}A(G)$. By Gilbert's Theorem, there are a Hilbert space $\Hilbert$ and bounded continuous $\boldsymbol{\xi}, \boldsymbol{\eta} \!: G \to \Hilbert$ such that (\ref{gileq}) holds and $\| \phi \|_{\Mcb} = \| \boldsymbol{\xi} \|_\infty \| \boldsymbol{\eta} \|_\infty$. Since
\[
  \check{\phi}(xy^{-1}) = \phi(y x^{-1} ) = \langle \boldsymbol{\xi}(y), \boldsymbol{\eta}(x) \rangle_\Hilbert
  = \overline{\langle \boldsymbol{\eta}(x),\boldsymbol{\xi}(y) \rangle_\Hilbert} 
  = \langle \boldsymbol{\eta}(x),\boldsymbol{\xi}(y) \rangle_{\overline{\Hilbert}} \qquad (x,y \in G),
\]
where $\overline{\Hilbert}$ denotes the complex conjugate Hilbert space of $\Hilbert$, it follows from Gilbert's Theorem that $\check{\phi} \in M_{cb}A(G)$ with $\| \check{\phi} \|_{\Mcb} \leq \| \boldsymbol{\xi} \|_\infty \| \boldsymbol{\eta} \|_\infty = \| \phi \|_{\Mcb}$.
\end{proof}
\par 
With Lemma \ref{checklem} at hand, we can prove a $A_{\Mcb}(G)$ version of \cite[Lemma 3.1]{Run}:
\begin{proposition} \label{amprop}
Let $G$ be a locally compact group such that $A_{\Mcb}(G)$ is $C$-amenable with $C \geq 1$. Then $\chi_\Gamma$ belongs to $M_{cb}A(G_d \times G_d)$ with $\| \chi_\Gamma \|_{\Mcb} \leq C$.
\end{proposition}
\begin{proof}
Let $( \boldsymbol{d}_\alpha )_\alpha$ be an approximate diagonal for $A_{\Mcb}(G)$ bounded by $C$. By Lemma \ref{checklem}, the net $( (\id \tensor \mbox{\ }^\vee)( (\boldsymbol{d}_\alpha) )_\alpha$ lies in $A_{\Mcb}(G) \tensor^\gamma A_{\Mcb}(G)$ and is also bounded by $C$. Obviously, $( (\id \tensor \mbox{\ }^\vee)( (\boldsymbol{d}_\alpha) )_\alpha$ converges to $\chi_\Gamma$ in the topology of pointwise convergence. Using more or less the same line of reasoning as in the proof of Theorem \ref{baithm}, we conclude that $\chi_\Gamma \in M_{cb}A(G_d \times G_d)$ with $\| \chi_\Gamma \|_{\Mcb} \leq C$.
\end{proof}
\begin{remark}
Let $A_M(G)$ be the closure of $A(G)$ in $\MA(G)$. The question for which $G$ the Banach algebra $A_M(G)$ is amenable seems to be more natural than the corresponding question for $A_{\Mcb}(G)$, but is apparently much less tractable (due to the fact that much less is known about $\MA(G)$ than about $M_{cb}A(G)$). For instance, we do not know whether or not an analog of Proposition \ref{amprop} holds for $A_M(G)$.
\end{remark}
\par 
Extending \cite[Theorem 3.5]{Run}, we obtain eventually:
\begin{theorem} \label{1am}
The following are equivalent for a locally compact group $G$:
\begin{items}
\item $G$ is abelian;
\item $A(G)$ is $1$-amenable;
\item $A_{\Mcb}(G)$ is $1$-amenable;
\item $A_{\Mcb}(G)$ is $C$-amenable with $C < \frac{2}{\sqrt{3}}$.
\end{items}
\end{theorem}
\begin{proof}
(i) $\Longleftrightarrow$ (ii) is \cite[Theorem 3.5]{Run} and (ii) $\Longrightarrow$ (iii) $\Longrightarrow$ (iv) are trivial.
\par 
(iv) $\Longrightarrow$ (i): If $A_{\Mcb}(G)$ is $C$-amenable with $C < \frac{2}{\sqrt{3}}$, then $\chi_\Gamma \in M_{cb}A(G_d \times G_d)$ is an idempotent with $\| \chi_\Gamma \|_{\Mcb} \leq C$ by Proposition \ref{amprop}. By Corollary \ref{idemcor}, this means that $\Gamma$ is a coset a of subgroup of $G \times G$ and thus a subgroup because $(e,e) \in \Gamma$. This is possible only if $G$ is abelian.
\end{proof}
\begin{remarks}
\item We do not know if the equivalent conditions in Theorem \ref{1am} are also equivalent to $A_M(G)$ being $1$-amenable.
\item In view of \cite[Theorem 2.3]{FR} and Theorem \ref{1am}, we believe that $A_{\Mcb}(G)$ is amenable if and only if $G$ is almost abelian. However, we have no proof in support of this belief. We do not even know whether or not $A_{\Mcb}(G)$ is amenable for $G = \free_2$.
\item As a consequence of Theorem \ref{1am}, we have for non-abelian $G$ that
\[
  \inf \{ C : \text{$A_{\Mcb}(G)$ is $C$-amenable} \} \geq \frac{2}{\sqrt{3}}
\]
(and possibly infinite). This, of course, entails that
\[
  \inf \{ C : \text{$A(G)$ is $C$-amenable} \} \geq \frac{2}{\sqrt{3}},
\]
which answers the question raised in the final remark of \cite{Run}.
\end{remarks}
\par 
We conclude the paper with an observation on amenable closed ideals of $A_{\Mcb}(G)$:
\begin{corollary} \label{amcor}
Let $G$ be a locally compact group, let $C \in \left[ 1, \frac{2}{\sqrt{3}} \right)$, and let $I$ be a non-zero, $C$-amenable, closed ideal of $A_{\Mcb}(G)$. Then $I$ is of the form $I(G \setminus xH)$, where $x \in G$ and $H$ is an open, abelian subgroup of $G$.
\end{corollary} 
\begin{proof}
Let $I$ be a non-zero, $C$-amenable, closed ideal of $A_{\Mcb}(G)$. From (\ref{diag2}), it is immediate that $I$ has an approximate identity bounded by $C$, and thus is of the form $I(G \setminus xH)$ for some open subgroup $H$ of $G$. In view of Proposition \ref{isoprop} and Theorem \ref{1am}, $H$ has to be abelian.
\end{proof}
\begin{remarks}
\item The restriction on $C$ in Corollary \ref{amcor} cannot be dropped: by \cite[(13) Bemerkung]{Lei}, there are infinite subsets $E$ of $\free_{2}$ such that $\chi_E M_{cb}A(G) \cong \ell^\infty(E)$, where $\cong$ stands for a not necessarily isometric isomorphism of Banach algebras. As $A_{\Mcb}(G)$ is Tauberian, it is then easy to see that the ideal $I = \chi_E A_{\Mcb}(G) = I(G \setminus E)$ is isomorphic to the commutative commutative $C^*$-algebra $c_0(E)$ and thus an amenable Banach algebra. Clearly, $I$ is not of the form described in Corollary \ref{amcor}. (It can be shown that $I$ is $4$-amenable and has an approximate identity bounded by $2$; see \cite{BFZ}.)
\item It is immediate from Corollary \ref{amcor} that $A_{\Mcb}(G)$ can have a non-zero, $C$-amenable, closed ideal if and only if $G$ has an open, abelian subgroup. In particular, for connected $G$, such ideals exist only if $G$ is abelian.
\end{remarks} 
\subsubsection*{Addendum} After this paper had been submitted we were informed by Ana-Maria Stan that Theorem 1.1 had been obtained independently in
\begin{quote}
\textsc{A.-M.\ Popa}, On idempotents of completely bounded multipliers of Fourier algebra $A(G)$. \textit{Indiana Univ.\ Math.\ J.}\ (to appear).
\end{quote}
\renewcommand{\baselinestretch}{1.0}
\dated
\vfill
\renewcommand{\baselinestretch}{1.2}
\begin{tabbing}
\textit{Second author's address}: \= Department of Mathematical and Statistical Sciences \kill 
\textit{First author's address}:  \> Department of Pure Mathematics \\
                                  \> University of Waterloo \\
                                  \> Waterloo, Ontario \\
                                  \> Canada N2L 3G1 \\[\medskipamount]
\textit{E-mail}:                  \> \texttt{beforres@math.ualberta.ca}\\[\bigskipamount]
\textit{Second author's address}: \> Department of Mathematical and Statistical Sciences \\
                                  \> University of Alberta \\
                                  \> Edmonton, Alberta \\
                                  \> Canada T6G 2G1 \\[\medskipamount]
\textit{E-mail}:                  \> \texttt{vrunde@ualberta.ca}\\[\medskipamount]
\textit{URL}:                     \> \texttt{http://www.math.ualberta.ca/$^\sim$runde/}
\end{tabbing}
\end{document}